\documentclass[12pt]{article}
\usepackage{latexsym, amscd, graphicx, amssymb}
\DeclareGraphicsExtensions{.ps}
    \title{{\bf  A functional-analytic theory of vertex  (operator)
algebras, II}}
    \author{Yi-Zhi Huang}
    \date{}
    
\begin{document}
    \bibliographystyle{alpha}
    \maketitle

   \newtheorem{thm}{Theorem}[section]
\newtheorem{defn}[thm]{Definition}
\newtheorem{prop}[thm]{Proposition}
\newtheorem{cor}[thm]{Corollary}
\newtheorem{rema}[thm]{Remark}
\newtheorem{lemma}[thm]{Lemma}
\newtheorem{app}[thm]{Application}
\newtheorem{prob}[thm]{Problem}
\newtheorem{conv}[thm]{Convention}
\newtheorem{conj}[thm]{Conjecture}
\newcommand{\halmos}{\rule{1ex}{1.4ex}}
\newcommand{\pfbox}{\hspace*{\fill}\mbox{$\halmos$}}

	\newcommand{\nno}{\nonumber}
	\newcommand{\lbar}{\bigg\vert}
\newcommand{\mbar}{\mbox{\large $\vert$}}
	\newcommand{\p}{\partial}
	\newcommand{\dps}{\displaystyle}
	\newcommand{\bra}{\langle}
	\newcommand{\ket}{\rangle}
 \newcommand{\res}{\mbox{\rm Res}}
\renewcommand{\hom}{\mbox{\rm Hom}}
\newcommand{\hol}{\mbox{\rm Hol}}
\newcommand{\dt}{\mbox{\rm Det}}
\newcommand{\edo}{\mbox{\rm End}\;}
 \newcommand{\pf}{{\it Proof.}\hspace{2ex}}
 \newcommand{\epf}{\hspace*{\fill}\mbox{$\halmos$}}
 \newcommand{\epfv}{\hspace*{\fill}\mbox{$\halmos$}\vspace{1em}}
 \newcommand{\epfe}{\hspace{2em}\halmos}
\newcommand{\nord}{\mbox{\scriptsize ${\circ\atop\circ}$}}
\newcommand{\wt}{\mbox{\rm wt}\ }
\newcommand{\swt}{\mbox{\rm {\scriptsize wt}}\ }
\newcommand{\clr}{\mbox{\rm clr}\ }

\begin{abstract}
For a finitely-generated 
vertex operator algebra $V$ of central charge $c\in \mathbb{C}$,
a locally convex topological 
completion $H^{V}$ is constructed. We construct on 
$H^{V}$ a structure of an
algebra over the operad of the $\frac{c}{2}$-th power $\dt^{c/2}$ 
of the determinant line 
bundle $\dt$ over  the moduli space of
genus-zero Riemann surfaces with ordered 
analytically parametrized boundary components. 
In particular, $H^{V}$ is a module for
the semi-group of the $\frac{c}{2}$-th power $\dt^{c/2}(1)$  of the
 determinant line bundle  over the moduli space of 
conformal equivalence classes
of annuli with analytically parametrized 
boundary components. 
The results in Part I for $\mathbb{Z}$-graded vertex algebras
are also reformulated in terms of the framed little disk operad.
Using May's recognition principle for double loop spaces,
one immediate consequence of such operadic formulations is 
that the  compactly generated spaces corresponding to (or 
the $k$-ifications of)
the locally convex completions constructed in Part I and in 
the present paper have the weak homotopy types of double loop 
spaces.
We also generalize the results above to locally-grading-restricted
conformal vertex algebras and to modules.
\end{abstract}

\renewcommand{\theequation}{\thesection.\arabic{equation}}
\renewcommand{\thethm}{\thesection.\arabic{thm}}
\setcounter{equation}{0}
\setcounter{thm}{0}
\setcounter{section}{-1}
\section{Introduction}

The present paper develops the functional-analytic aspects of
vertex operator algebras. More specifically, we construct
a locally convex topological completion of a finitely-generated
vertex operator algebra  and
a structure on this completion of an algebra over a certain natural operad
constructed from
genus-zero Riemann surfaces with boundaries.
We obtain representation-theoretic and homotopy-theoretic
consequences 
and give generalizations to more general algebras and modules.

For a complex number $c$, consider the
sequence $\dt^{c/2}$ of the $\frac{c}{2}$-th power of the 
determinant line bundles $\dt^{c/2}(n)$, $n\ge 0$, over 
the moduli spaces of
genus-zero Riemann surfaces with $n+1$ ordered 
analytically parametrized boundaries. This sequence $\dt^{c/2}$
has a natural structure of (genuine) operad. (See \cite{M}, \cite{HL1}, 
\cite{HL2} and
Appendix C of \cite{H3} for the notion of operads and other 
related notions and see \cite{S1}, \cite{S2} 
and Appendix D of \cite{H3} for determinant line bundles.)
An algebra 
over $\dt^{c/2}$ such that the underlying vector space 
is a complete locally convex topological vector space and the 
corresponding maps are continuous and depend holomorphic on 
$\dt^{c/2}$
is called a {\it genus-zero holomorphic conformal field theory
of central charge $c$}. See \cite{S1} and \cite{S2} for 
a geometric definition of conformal field theory 
in the more general case of arbitrary genus and nonholomorphic theories.

Genus-zero conformal field theories are the starting point of 
a number of papers on algebraic structures derived from
conformal field theories (see, for example, 
\cite{KSV} \cite{KVZ}). 
But the construction of examples of conformal field theories, 
even in this genus-zero case, 
is difficult and subtle. It has been expected that vertex operator 
algebras will give examples of such genus-zero holomorphic theories. But 
it is clear that vertex operator algebras themselves are 
not such theories. In fact, in \cite{H1}, \cite{H2}
and \cite{H3}, it was established 
that a vertex operator algebra has only the structure of an algebra over a 
$\mathbb{C}^{\times}$-rescalable partial operad 
in the sense of \cite{HL1} and \cite{HL2}.
It is also clear that to construct such a theory from a vertex operator
algebra, one first has to construct a suitable locally convex completion
of the algebra.
We know that $\dt^{c/2}$ is generated by $\dt^{c/2}(1)$ and $\dt^{c/2}(2)$,
the $\frac{c}{2}$-th power of the  determinant line
bundle over of genus-zero Riemann surfaces with two and
three, respectively, ordered analytically 
parametrized boundary components. Thus one must next 
construct continuous linear maps associated to 
elements in $\dt^{c/2}(1)$ and $\dt^{c/2}(2)$. Combining these maps 
with the geometric formulation of 
vertex operator algebras in terms of partial operads in \cite{H3}, 
it is easy to see that we will have a genus-zero holomorphic
conformal field theory.

The main purpose of the present paper is to carry out this
construction of  genus-zero holomorphic 
conformal field theories from finitely-generated vertex operator
algebras. The results in Part I for 
finitely-generated $\mathbb{Z}$-graded 
vertex algebras
are also reformulated in terms of the framed little disk operad.
Note that any genus-zero conformal field theory must be a
representation of the semi-group $\dt^{c/2}(1)$,
the $\frac{c}{2}$-th power of the determinant
line bundle over the moduli space of annuli with analytically
parametrized boundary components. Thus, in particular, we 
construct in this paper 
a representation of $\dt^{c/2}(1)$ from a finitely-generated
vertex operator algebra. 
In fact, from the construction
it is easy to see that part of our construction 
actually gives a representation of $\dt^{c/2}(1)$
from an arbitrary $\mathbb{Z}$-graded representation of 
the Virasoro algebra satisfying a certain truncation condition. 
As far as the author knows, there seems to be  no such general results
on the integration of representations of the Virasoro algebra in 
the literature.

Combining the operadic formulations mentioned above with May's
recognition principle for double loop spaces \cite{M}, we conclude
that the compactly generated spaces corresponding to (or the
$k$-ifications of) the locally convex completions constructed in Part
I and in the present paper have the weak homotopy types of double loop
spaces. It is known that vertex operator algebras are a basic
ingredient in conformal field theory and that conformal field theories
describe string theory or M theory perturbatively.  In string theory,
there are two kinds of geometry involved, the "world-sheet" geometry
and "space-time" geometry.  The operad $\dt^{c/2}$ is part of the
world-sheet geometry. The double loop space structures are interesting
because they give us some "space-time" information about the vertex
(operator) algebra. Since the operad $\dt^{c/2}$ has a much richer
structure than the little disk operad, one should be able to recognize
more properties of algebras over it. It will be much more interesting
if one can recognize topological properties homeomorphically, not just
(weak) homotopically, or even recognize some geometric properties. It
will be especially interesting to see what geometric and topological
properties can be recognized from structures 
associated to conformal field theories such as the minimal models
which are constructed without any "space-time" geometry information.

These constructions and results above generalize 
to locally-grading-res\-tricted
conformal vertex algebras without any difficulty.
These generalizations have been used in \cite{HZ}.
We also give the corresponding  results
for modules without giving 
detailed proofs.

The present paper is organized as follows: 
In Section 1, a locally convex topological completion $H^{V}$
of a finitely-generated vertex 
operator algebra $V$ of central charge $c\in \mathbb{C}$
is constructed. In Section 2, 
a structure of a representation on $H^{V}$ of 
$\dt^{c/2}(1)$ is constructed. In Section 3, 
we construct linear continuous maps from
the completed tensor product of two copies of 
$H^{V}$ to $H^{V}$ associated 
to elements of $\dt^{c/2}(2)$. In Section 4, we first reformulate 
the result in Part I (\cite{H4}) in terms of the framed little disk 
operad. Then
we state the main result (Theorem \ref{main}) of the present paper.
Structures of  double loop spaces on the 
compactly generated spaces corresponding to (or the $k$-ifications of)
the completions constructed in \cite{H4}
and in this paper are also stated in this section.
The statements of the generalizations to locally-grading-restricted
conformal vertex algebras and the corresponding results 
for modules are given in 
Sections 5 and 6, respectively.

\paragraph{Acknowledgment} I am grateful to J. Peter May and 
Nick Kuhn
for discussions on compactly generated spaces and 
the recognition principle for 
double loop spaces. This research is supported in part 
by NSF grant DMS-0070800. 

\renewcommand{\theequation}{\thesection.\arabic{equation}}
\renewcommand{\thethm}{\thesection.\arabic{thm}}
\setcounter{equation}{0}
\setcounter{thm}{0}

\section{A locally convex completion of a finitely-generated vertex 
operator algebra}

In this section, we construct a locally convex topological completion
of a finitely-generated 
vertex operator algebra $V$. The topological completion is larger
than the topological completion of a finitely-generated 
$\mathbb{Z}$-graded
grading-restricted vertex algebra constructed in Part I (\cite{H4}). For
simplicity, we shall use the same notation $H^{V}$ as in \cite{H4} to
denote the topological completion we shall construct in the present
paper. But we warn the reader that $H^{V}$ in the present paper is
larger than $H^{V}$ in \cite{H4}. As in \cite{H4}, since $V$ is 
fixed in the present paper, we shall denote $H^{V}$ simply by $H$.

First we need to consider some geometric objects. A {\it disk} is 
a genus-zero Riemann surface with a connected boundary.
A smooth invertible map 
from $S^{1}$ to the boundary of a disk is called an {\it analytic 
parametrization} if it can be extended to an analytic map from a 
neighborhood of $S^{1}$ inside  the closed unit disk on the complex plane
to a neighborhood of the boundary of the disk.
A {\it disk with analytically parametrized boundary} is 
a disk equipped with an analytic parametrization of its boundary.
For $k\ge 0$, a {\it $k$-punctured 
disk with analytically parametrized boundary} 
is a disk with analytically parametrized boundary
and  $k$ ordered and distinct
 points in the interior of the disk. {\it Conformal
equivalences} between $k$-punctured 
disks with analytically parametrized boundaries are defined in the 
obvious way.

Let $\Theta(k)$, $k\ge 0$, be the moduli spaces 
of $k$-punctured disks with
analytically parametrized boundaries and let $\Theta=\cup_{k\ge
0}\Theta(k)$. Also consider the moduli spaces $\mathcal{B}_{0, 1, k}$, $k\ge
0$, of genus-zero Riemann surfaces with ordered analytically parametrized
boundary components, one positively oriented and the other negatively
oriented and ordered. The sequence $\{\mathcal{B}_{0, 1, k}\}_{k\ge 0}$
has a natural structure of an analytic operad and this operad is
isomorphic to the suboperad $K_{\mathfrak{H}_{1}}$ of the sphere partial
operad $K$ discussed in Section 6.4 of \cite{H3}. It is clear that $\Theta$
has a natural structure of a space over the operad $\{\mathcal{B}_{0, 1,
k}\}_{k\ge 0}$ or equivalently over the operad $K_{\mathfrak{H}_{1}}$.

For any $k\ge 0$, we have an injective map 
from $\Theta(k)$ to $K(k)$ defined as
follows: Take any element of $\Theta(k)$, that is, a conformal equivalence
class of $k$-punctured disks with analytically parametrized
boundaries. For any $k$-punctured disk with analytically parametrized
boundary in this conformal equivalence class, by sewing the 
union of the exterior of $S^{1}$ and $\infty$
to this $k$-punctured disk using the analytic boundary
parametrization, we obtain a $k+1$-punctured genus-zero Riemann
surface, one puncture negatively oriented and the other puncture
positively oriented and ordered, together with a local analytic
coordinate vanishing at the negatively oriented puncture. Using the
uniformization theorem, this $k+1$-punctured genus-zero Riemann
surface with a local coordinate at the negatively oriented puncture
is conformally equivalent to $\mathbb{C}\cup \{\infty\}$ with $k+1$
punctures, one negatively oriented and the other positively oriented
and ordered, together with a local coordinate vanishing at the
negatively oriented puncture. Moreover, we can choose the conformal
equivalence (analytic diffeomorphism) such that the negatively
oriented puncture is mapped to $\infty$, the $k$-th positively
oriented puncture is mapped to $0$ and the derivative at $\infty$ of
the local coordinate map vanishing at $\infty$ is $1$. Adding the
standard local coordinates vanishing at the positively oriented
punctures, we obtain a canonical sphere with tubes of type $(1, k)$
(see Chapter 3 of \cite{H3}).  It is clear that this canonical sphere
with tubes of type $(1, k)$ is independent of the choice of the
$k$-punctured disk with analytically parametrized boundary in the given
conformal equivalence class.  We define the image of the element of
$\Theta(k)$ to be the conformal equivalence class of spheres with tubes of
type $(1, k)$ containing this canonical sphere with tubes of type $(1,
k)$. We obtain a map from $\Theta(k)$ to $K(k)$. Clearly this map is
injective. We shall identify $\Theta(k)$ with its image in $K(k)$. Note
that any element of $\Theta(k)$ viewed as an element of $K(k)$ is of the form 
\begin{equation}\label{xn-elt}
P=(z_{1}, \dots, z_{k-1}; A, (1, \mathbf{0}), \dots, (1, \mathbf{0}))
\end{equation}
where $A\in H$. But also note that not all elements of $K(k)$ of this form
is an element of $\Theta(k)$.

Note that $\Theta(k)$ can also be viewed as a subset of the Banach space
$\mathbb{C}^{k-1}\times \hol(D^{1})$ where $\hol(D^{1})$ is
the Banach space of all functions continuous on the closed unit disk $D^{1}$
and holomorphic on the open unit disk. We give $\Theta(k)$
the topology and analytic structure induced from 
those on $\mathbb{C}^{k-1}\times \hol(D^{1})$.

Let $(V, Y, \mathbf{1}, \omega)$ be a vertex operator algebra (in the
sense of \cite{FLM} and \cite{FHL}).  By the isomorphism theorem
proved in Chapter 5 of \cite{H3}, there exists a canonical geometric
vertex operator algebra structure on $V$. Let $\nu_{k}: K(k)\to
\hom(V^{\otimes k}, \overline{V})$, $k\ge 0$,  be the maps defining 
the geometric
vertex operator algebra structure on $V$. Then  for
any $v'\in V'$, any $u_{1}, \dots, u_{k}, v\in V$,
$$\langle v', (\nu_{k}(P))(u_{1}\otimes \cdots \otimes u_{k}\otimes 
v )\rangle$$
as a function of $P$ is 
meromorphic  on $K(k)$. 
Thus for any $u_{1}, \dots, u_{k}, v\in V$ and any 
$P\in K(k)$, we have an element 
$$Q(u_{1}, \dots, u_{k}, v; P)
=(\nu_{k}(P))(u_{1}\otimes \cdots \otimes u_{k}\otimes 
v )\in \overline{V}.$$
In particular, for any $u_{1}, \dots, u_{k}, v\in V$ and any 
$P\in \Theta(k)$, we have an element 
$Q(u_{1}, \dots, u_{k}, v; P)\in \overline{V}$ since 
$\Theta(k)$ can be viewed as a subset of $K(k)$.

For $k\ge 0$ and $n>0$, let 
\begin{eqnarray*}
J_{n}^{(k)}&=&\{(z_{1}, \dots, z_{k-1}; A, (1, \mathbf{0}), 
\dots, (1, \mathbf{0}))
\in \Theta(k)\;|\; |z_{i}-z_{j}|\ge \frac{1}{n}, i\ne j, \nno\\
&&\quad |z_{i}|>\frac{1}{n}, i=1, 
\dots, k,
\mbox{\rm the distances from 
$z_{i}$, $i=1, \dots, k-1$,} \nno\\
&&\quad \mbox{\rm $0$ 
to $C_{1}$ and from $0$ to $C_{1}^{-1}$ 
are large than or equal to ${\displaystyle \frac{1}{n}}$}\},
\end{eqnarray*}
where 
$$C_{1}=f_{A}(\{w\in \mathbb{C}\;|\;|w|=1\},$$
$$f_{A}(w)=e^{
\sum_{j>0}A_{j}w^{j+1}\frac{d}{dw}} w$$
and
$$C^{-1}_{1}=\{w^{-1}\;|\; w\in C_{1}\}.$$
Then we see 
$\Theta(k)=\cup_{n>0}J_{n}^{(k)},\;\;\;k\ge 0.$

We denote the projections from $V$ to $V_{(n)}$, $n\in
\mathbb{Z}$, by $P_{n}$ as in \cite{H3}.
For fixed $k\ge 0$, by the sewing axiom for geometric vertex operator algebras
in \cite{H3}, 
\begin{equation}\label{conv-fk1}
\sum_{n\in \mathbb{Z}}\langle v', (\nu_{l}(Q))(v_{1}\otimes \cdots
\otimes v_{l-1}\otimes
(P_{n}(Q(u_{1}, \dots, u_{k}, v; P)))\rangle,
\end{equation}
$v'\in V'$, $u_{1}, \dots, u_{k}, v_{1}, \dots, v_{l}\in V$,
$P\in \Theta(k)$ and $Q\in K_{\mathfrak{H}_{1}}(l)$,
is absolutely convergent. 
For fixed $v'\in V'$, $u_{1}, \dots, u_{k}, v_{1}, \dots, v_{l}\in V$,
and $Q\in K_{\mathfrak{H}_{1}}(l)$, the sum of (\ref{conv-fk1}) give
a function on $\Theta_{k}$. 

\begin{lemma}
The functions defined by  the sum of (\ref{conv-fk1})
is bounded 
on $J_{n}^{(k)}$, $n>0$. 
\end{lemma}
\pf
By the sewing axiom for geometric vertex operator algebras
in \cite{H3},
(\ref{conv-fk1}) is equal to 
\begin{equation}\label{conv-fk1.5}
\langle v', (\nu_{k+l-1}(Q_{^{l}}\infty_{^{0}}P))(v_{1}\otimes 
\cdots \otimes v_{l-1}\otimes u_{1}\otimes \cdots\otimes u_{k}\otimes v)
\rangle.
\end{equation}
From (\ref{conv-fk1.5}) and 
the definition of $\nu_{k+l-1}$ in \cite{H3}, we see that
to prove the lemma, we need only show that when $P\in J_{n}^{(k)}$, 
the distances between distinct punctures of $Q_{^{l}}\infty_{^{0}}P$ are
larger than  a fixed positive number depending only on $n$, and 
each expansion coefficient, as a function of $P$, of the analytic 
local coordinate maps vanishing at these punctures are 
bounded on $J_{n}^{(k)}$.

We first recall some facts and results from \cite{H3}.
Let 
$$Q=(\xi_{1}, \dots, \xi_{l-1}; B^{(0)}, (b^{(1)}_{0}, B^{(1)}),
\dots, (b_{0}^{(l)}, B^{(l)}))$$
and 
$$f_{B^{(i)}, b_{0}^{(i)}}(w)=b_{0}^{(i)}e^{
\sum_{j>0}B^{(i)}_{j}w^{j+1}\frac{d}{dw}} w$$
for $i=0, \dots, l$. 
We shall also use the same notations
$f_{B^{(i)}, b_{0}^{(i)}}$, $i=0, \dots, l$, and $f_{A}$ to
denote the corresponding local coordinate maps. 
Then by the study of the sewing operation in 
\cite{H3}, the sewing equation 
$$F^{(1)}(w)=F^{(2)}\left(f_{A}^{-1}\left(\frac{1}
{f_{B^{(l)}, b_{0}^{(l)}}(w)}\right)\right)$$
together with the normalization conditions
\begin{eqnarray*}
F^{(1)}(\infty)&=&\infty,\\
F^{(2)}(0)&=&0,\\
\lim_{w\to \infty}\frac{F^{(1)}}{w}&=&1
\end{eqnarray*}
has a unique solution pair $F^{(1)}, F^{(2)}$ and 
the positively oriented punctures of $Q_{^{l}}\infty_{^{0}}P$ 
corresponding to the positively oriented punctures of $P$ are 
$F^{(2)}(z_{1})$, $\dots, F^{(2)}(z_{k-1})$ and $0$.
The local coordinate maps vanishing at these punctures are 
$F^{(2)}(w)-F^{(2)}(z_{1}), \dots, F^{(2)}(w)-F^{(2)}(z_{l})$
and $F^{(2)}(w)$, respectively. 
It is also proved in \cite{H3} that the sewing operation 
is analytic. In particular, it is continuous. Thus 
$Q_{^{l}}\infty_{^{0}}P$ is continuous in $P\in \Theta(k)$.
In fact the proof actually proves that $F^{(1)}$ and $F^{(1)}$ 
depend on $f_{A}$ and $f_{B^{(l)}, b_{0}^{(l)}}$ 
analytically and in particular 
continuously.

First we prove that  when $P\in J_{n}^{(k)}$, the distances between 
distinct punctures of $Q_{^{l}}\infty_{^{0}}P$ are
larger than  a fixed positive number depending only on $n$.
If this is not true, then there is a sequence $\{P_{m}\}_{m>0}$
in $J_{n}^{(k)}$ and two punctures on 
$Q_{^{l}}\infty_{^{0}}P_{m}$ for each $m>0$ having the 
same orders, such that the distance between 
these two punctures goes to $0$ when $m$ goes to $\infty$. 
We consider the case that these two punctures are positively oriented 
punctures corresponding to two nonzero positively oriented punctures
on $P_{m}$. If we use $z_{1}(P_{m}), \dots, z_{k-1}(P_{m})$
to denote the punctures of $P_{m}$ and $F_{m}^{(2)}$ and $F_{m}^{(2)}$
the  solution of the sewing equation and the normalization conditions
with $P$ replaced by $P_{m}$, then by the results in \cite{H3}
we recalled above, 
these two punctures $Q_{^{l}}\infty_{^{0}}P_{m}$
must be of the form $F_{m}^{(2)}(z_{p}(P_{m}))$ and 
$F_{m}^{(2)}(z_{q}(P_{m}))$
for some $0<p, q<k$. 

On the other hand,
we can also obtain  $z_{1}(P_{m}), \dots, z_{k-1}(P_{m})$ from 
$$F_{m}^{(2)}(z_{1}(P_{m})), \dots, F_{m}^{(2)}(z_{k-1}(P_{m}))$$
as follows: 
We sew the 
first puncture of 
$(\mathbf{0}, (b_{0}^{(l)}, B^{(l)}(b_{0}^{(1)})))$ to 
the $0$-the puncture of 
$$(F_{m}^{(2)}(z_{1}(P_{m})), \dots, F_{m}^{(2)}(z_{k-1}(P_{m})); -\Psi^{-}, 
(1, \mathbf{0}), \dots, (1, \mathbf{0})),$$
where 
$$B^{(l)}(b_{0}^{(l)})=\{(b_{0}^{(l)})^{j}B^{(l)}_{j}\}_{j>0}$$
and $-\Psi^{-}=\{-\Psi_{j}\}_{j<0}$ is the sequence defined by 
$$F^{(1)}(w)
=e^{-\sum_{j<0}\Psi_{j}w^{j+1}\frac{d}{dw}}w.$$
Then the positively oriented punctures of the resulting 
element is $z_{1}(P_{m}), \dots$, $z_{k}(P_{m})$. In particular, 
$z_{1}(P_{m}), \dots, z_{k}(P_{m})$ depend continuously
on 
$$F_{m}^{(2)}(z_{1}(P_{m})), 
\dots, F_{m}^{(2)}(z_{k-1}(P_{m})).$$
Thus since the distance between $F_{m}^{(2)}(z_{p}(P_{m}))$ and 
$F_{m}^{(2)}(z_{q}(P_{m}))$ goes to $0$ when $m$ goes to $\infty$,
the distance between $z_{p}(P_{m}$ and $z_{q}(P_{m}$ must also 
goes to $0$ when $m$ goes to $\infty$. But $\{P_{m}\}_{m>0}$ is 
in $J_{n}^{(k)}$ and by the definition of $J_{n}^{(k)}$, this is 
impossible. Similarly we get contradictions in the other cases. 
Thus when $P\in J_{n}^{(k)}$, the distances between 
distinct punctures of $Q_{^{l}}\infty_{^{0}}P$ are
larger than  a fixed positive number depending only on $n$.

Now we prove that each expansion coefficient, as
a function of $P$,  of the analytic 
local coordinate maps vanishing at the punctures of $Q_{^{l}}\infty_{^{0}}P$
is bounded on $J_{n}^{(k)}$. For simplicity, we prove this for the 
expansion coefficients of the analytic 
local coordinate map vanishing at the last puncture $0$ of 
$Q_{^{l}}\infty_{^{0}}P$. By the results in \cite{H3} we recalled 
above, the local coordinate map vanishing at $0$ is $(F^{(2)})^{-1}$.
Since the expansion coefficients of $(F^{(2)})^{-1}$ at $0$ are
polynomials in the expansion coefficients of $F^{(2)}$, we
need only show that each expansion coefficient of $F^{(2)}$ 
as a function of $P$ is bounded on $J_{n}^{(k)}$. Note that 
the domain of $F^{(2)}$ contains $C_{1}$ and 
the interior of $C_{1}$. When 
$P\in J_{n}^{(k)}$, the union of $C_{1}$ and 
the interior of $C_{1}$ always contains
the closed disk centered at $0$ of radius $1/n$. 
Since the radius $1/n$ is independent of $P$, 
using the Cauchy formulas for the expansion coefficients
of $F^{(2)}$, we see that each of these coefficients is
bounded on $J_{n}^{(k)}$. 
\epfv

Let $\tilde{G}$ be the subspace of $V^{*}$ 
consisting of linear functionals $\lambda$ on $V$
such that for any $k\ge 0$,
$u_{1}, 
\dots, u_{k}, v\in V$, $P\in \Theta(k)$,
\begin{equation}\label{conv-fk2}
\sum_{n\in \mathbb{Z}}\lambda(P_{n}(Q(u_{1}, \dots, u_{k}, v; P)))
\end{equation}
is absolutely convergent  and its
sum  as a function on $\Theta(k)$ is bounded on $J_{n}^{(k)}$, 
$n>0$. 
The dual pair $(V^{*},
V)$ of vector spaces gives $V^{*}$ a locally convex topology.  With
the topology induced from the one on $V^{*}$, $\tilde{G}$ is also a
locally convex space. Note that $V'$ is a subspace of $\tilde{G}$. (We
warn the reader that $\tilde{G}$ here is different from $\tilde{G}$ in
\cite{H4}. In the present paper, many notations we use are the same as
the corresponding notations in \cite{H4}. But what they denote 
are different from what the same notations denote in \cite{H4}.)

We denote the 
analytic function on $\Theta(k)$ defined by (\ref{conv-fk2})
by
$$g_{k}(\lambda \otimes u_{1}\otimes \cdots\otimes u_{k}
\otimes v)$$
since (\ref{conv-fk2}) is multilinear in $\lambda$, $u_{1}, \dots, u_{k}$
and $v$. 
These functions  span a vector space $F_{k}$ of analytic functions on 
$\Theta(k)$. We obtain a linear map 
$$g_{k}: \tilde{G}\otimes V^{\otimes (k+1)}\to F_{k}.$$
By definition,  elements of $F_{k}$ are bounded on $J_{n}^{(k)}$,
$n>0$.

We define a family of norms $\|\cdot\|_{F_{k}, n}$, $n>0$, 
on $F_{k}$
by 
$$\|g\|_{F_{k}, n}=\sup_{Q\in J^{(k)}_{n}}g(Q)$$
for $g\in F_{k}$. These norms 
give a locally convex topology on $F_{k}$. Note that a net 
$\{f_{\alpha}\}_{\alpha\in \mathcal{A}}$ (where $\mathcal{A}$ is 
an index set) in $F_{k}$ is convergent to $f\in \Theta_{k}$ if and only if
it is convergent uniformly in $J^{(k)}_{n}$ for $n>n_{0}$ where 
$n_{0}$ is a positive integer.

For any $k\ge 0$, there is  an embedding $\iota_{F_{k}}$ from
$F_{k}$ to $F_{k+1}$ defined as follows:  We use 
$$(z_{0}, \dots,
z_{k-1}; A, (1, \mathbf{0}), \dots, (1, \mathbf{0}))$$ instead of 
$$(z_{1}, \dots, z_{k}; A, (1, \mathbf{0}), \dots, (1, \mathbf{0}))$$
to denote the elements of
$\Theta_{k+1}$. 
For $\lambda\in \tilde{G}$, $u_{1}, \dots, u_{k}, v\in V$, since 
$$Y(\mathbf{1}, z)=1$$
for any nonzero complex number $z$,
$$g_{k+1}(\lambda \otimes  {\bf
1}\otimes u_{1}\otimes \cdots \otimes u_{k}\otimes 
v)$$ 
as a function of $(z_{0}, \dots, z_{k-1}; A, (1, \mathbf{0}), \dots,
(1, \mathbf{0}))$ 
is in fact independent of $z_{0}$, and is equal to
$$g_{k}(\lambda \otimes u_{1}\otimes \cdots \otimes u_{k}\otimes 
v)$$ 
as a function in $(z_{1}, \dots, z_{k-1}; A, (1, \mathbf{0}), \dots,
(1, \mathbf{0}))$. Thus we obtain a
well-defined linear
map 
$$\iota_{F_{k}}: F_{k}\to F_{k+1}$$
such that 
$$\iota_{F_{k}}\circ g_{k}=g_{k+1}\circ \phi_{k}$$
where 
$$\phi_{k}: \tilde{G}\otimes V^{\otimes (k+1)}  \to \tilde{G} \otimes 
V^{\otimes (k+2)}$$
is defined by
$$\phi_{k}(\lambda \otimes  u_{1}\otimes \cdots \otimes u_{k}\otimes v)
=\lambda\otimes \mathbf{1}\otimes 
  u_{1}\otimes \cdots \otimes u_{k}\otimes v$$
for $\lambda\in \tilde{G}$, $u_{1}, \dots, u_{k},
v\in V$. It is clear
that $\iota_{F_{k}}$ is injective. Thus we can regard $F_{k}$ as a
subspace of $F_{k+1}$. Moreover, we have:

\begin{prop}\label{k:k+1}
For any $k\ge 0$,   $\iota_{F_{k}}$ as a map from $F_{k}$ to 
$\iota_{F_{k}}(F_{k})$
is continuous and open. In other words, the topology on $F_{k}$
is induced from that on $F_{k+1}$.
\end{prop}
\pf
We consider the two topologies on $F_{k}$, one is the topology
defined above for $F_{k}$ 
and the other induced from the topology on $F_{k+1}$.
We need only prove that for any $n>0$, (i) the norm
$\|\cdot\|_{F_{k}, n}$ is continuous in the topology induced from the one
on $F_{k+1}$, and (ii) the restriction of the norm
$\|\cdot\|_{F_{k+1}, n}$ 
to $F_{k}$ is continuous in the topology on
$F_{k}$. 

Let 
$\{f_{\alpha}\}_{\alpha\in \mathcal{A}}$ (where $\mathcal{A}$ is 
an index set) be a net in $F_{k}$
convergent in the topology induced from the one on $F_{k+1}$.  Then
$\{f_{\alpha}\}_{\alpha\in \mathcal{A}}$, when viewed as a net of
functions in 
$$(z_{0}, z_{1}, \dots, z_{k-1}; A, (\mathbf{0}, (1, \mathbf{0}),
\dots, (\mathbf{0}, (1, \mathbf{0}))$$
is convergent uniformly
on $J^{(k+1)}_{n}$, $n>0$. Since $f_{\alpha}$, $\alpha\in \mathcal{A}$, 
are independent of $z_{0}$, $\{f_{\alpha}\}_{\alpha\in \mathcal{A}}$ 
is in fact convergent uniformly on the sets $J_{n+1}^{(k)}$, $n>0$,
proving (i). Now let $\{f_{\alpha}\}_{\alpha\in \mathcal{A}}$ 
be a net in $F_{k}$ convergent in the topology on $F_{k}$. Then
$\{f_{\alpha}\}_{\alpha\in \mathcal{A}}$ is convergent uniformly on 
$J_{n}^{(k)}$, $n>0$. If we view $f_{\alpha}$, $\alpha\in \mathcal{A}$,
as functions on $\mathbb{C}\times \Theta(k)$, then the net
$\{f_{\alpha}\}_{\alpha\in \mathcal{A}}$ is convergent uniformly on 
$(\mathbb{C}\times J_{n}^{(k)})\cap \Theta(k+1)$ (where we view 
$\Theta(k+1)$ as a subset of $\mathbb{C}\times \Theta(k)$.
Since $J_{n}^{(k+1)}\subset (\mathbb{C}\times J_{n}^{(k)})\cap \Theta(k+1)$,
$\{f_{\alpha}\}_{\alpha\in \mathcal{A}}$
is convergent uniformly on $J_{n}^{(k+1)}$, $n>0$,
proving (ii). 
\epfv

We equip the topological dual space $F_{k}^{*}$, $k\ge 0$, of $F_{k}$
with the strong topology, that is, the topology of uniform convergence 
on all the weakly bounded subsets of $F_{k}$.
Then $F_{k}^{*}$ is a locally convex space.

For  $k\ge 0$, we define a linear map 
$$\gamma_{k}: F_{k+1}\to F_{k}$$
as follows: We use $P=(z_{0}, \dots,
z_{k-1}; A, (1, \mathbf{0}), \dots, (1, \mathbf{0}))$
to denote an element of
$\Theta_{k+1}$.
Recall that 
$$C_{1}=f_{A}(\{w\in \mathbb{C}\;|\;|w|=1\}$$
and
$$f_{A}(w)=e^{\sum_{j>0}A_{j}w^{j+1}\frac{d}{dw}}w.$$
We define
\begin{eqnarray}\label{gamma-k-1}
\lefteqn{\gamma_{k}(g_{k+1}(\lambda\otimes 
u_{0}\otimes u_{1}\otimes \cdots \otimes 
u_{k}\otimes v))}\nno\\
&&=\frac{1}{2\pi \sqrt{-1}}\oint_{C_{1}}z_{0}^{-1}
g_{k+1}(\lambda\otimes 
u_{0}\otimes u_{1}\otimes \cdots \otimes 
u_{k}\otimes v)dz_{0}
\end{eqnarray}
for
$\lambda\in \tilde{G}$, $u_{0}, u_{1}, \dots, u_{k}, v\in V$. 

We still need to show that the right-hand side of (\ref{gamma-k-1})
 is indeed in $F_{k}$. Let $P'=(z_{1}, \dots,
z_{k-1}; A, (1, \mathbf{0}), \dots, (1, \mathbf{0}))\in \Theta_{k}$.
Then we have 
\begin{equation}\label{decomp}
P=(f_{A}^{-1}(z_{0}); \mathbf{0}, B(z_{0}), 
(1, \mathbf{0}))_{^{2}}\infty_{^{0}}P'
\end{equation}
(see formula (A.6.1) in \cite{H3}),
where 
$$B(z_{0})=\hat{E}^{-1}\left(\frac{1}{f_{A}^{-1}
\left(\frac{1}{x+\frac{1}{f_{A}(z_{0})}}\right)}-\frac{1}{z_{0}}\right).$$
By the definition of $g_{k}$ and (\ref{decomp}), we have
\begin{eqnarray}\label{gamma-k-2}
\lefteqn{\gamma_{k}(g_{k+1}(\lambda\otimes
u_{0}\otimes u_{1}\otimes \cdots \otimes 
u_{k}\otimes v))}\nno\\
&&=\frac{1}{2\pi \sqrt{-1}}\oint_{C_{1}}z_{0}^{-1}
\sum_{n\in \mathbb{Z}}\lambda(P_{n}(Q(u_{0}, \dots, u_{k}, v; P)))
dz_{0}.
\end{eqnarray}
Since the series 
$$\sum_{n\in \mathbb{Z}}\lambda(P_{n}(Q(u_{0}, \dots, u_{k}, v; P)))$$
is absolutely convergent,
the right-hand side of (\ref{gamma-k-2}) is equal to
\begin{eqnarray}\label{gamma-k-3}
\lefteqn{\frac{1}{2\pi \sqrt{-1}}\sum_{n\in \mathbb{Z}}
\oint_{C_{1}}z_{0}^{-1}
\lambda(
P_{n}(\nu_{k+1}(P)(u_{0}\otimes
u_{1}\otimes \cdots \otimes u_{k}\otimes 
v )))dz_{0}}\nno\\
&&=\frac{1}{2\pi \sqrt{-1}}\sum_{n\in \mathbb{Z}}\oint_{C_{1}}z_{0}^{-1}
\lambda(P_{n}(
(\nu_{k+1}((f^{-1}_{A}(z_{0}); \mathbf{0}, B(z_{0}), 
(1, \mathbf{0}))_{^{2}}\infty_{^{0}}P')\nno\\
&&\hspace{6em}(u_{0}\otimes
u_{1}\otimes \cdots \otimes u_{k}\otimes 
v )))dz_{0}\nno\\
&&=\frac{1}{2\pi \sqrt{-1}}\sum_{n\in \mathbb{Z}}\oint_{C_{1}}z_{0}^{-1}
\lambda(P_{n}(
(\nu_{2}((f^{-1}_{A}(z_{0}); \mathbf{0}, B(z_{0}), 
(1, \mathbf{0})))_{^{2}}*_{^{0}}\nu_{k}(P'))\nno\\
&&\hspace{6em}(u_{0}\otimes
u_{1}\otimes \cdots \otimes u_{k}\otimes 
v )))dz_{0}\nno\\
&&=\frac{1}{2\pi \sqrt{-1}}\sum_{n\in \mathbb{Z}}\oint_{C_{1}}z_{0}^{-1}
\lambda(Y(e^{-\sum_{j>0}B(z_{0})L(j)}
u_{0}, f^{-1}_{A}(z_{0}))\nno\\
&&\hspace{6em}P_{n}((\nu_{k}(P'))(
u_{1}\otimes \cdots \otimes u_{k}\otimes 
v )))dz_{0}\nno\\
&&=\frac{1}{2\pi \sqrt{-1}}\sum_{n\in \mathbb{Z}}
\oint_{|w|=1}(f_{A}(w))^{-1}f'_{A}(w)
\lambda(Y(e^{-\sum_{j>0}B(f_{A}(w))L(j)}
u_{0}, w)\nno\\
&&\hspace{6em}P_{n}((\nu_{k}(P'))(
u_{1}\otimes \cdots \otimes u_{k}\otimes 
v )))dw\nno\\
&&=\sum_{n\in \mathbb{Z}}
\res_{w}(f_{A}(w))^{-1}f'_{A}(w)
\lambda(Y(e^{-\sum_{j>0}B(f_{A}(w))L(j)}
u_{0}, w)\nno\\
&&\hspace{6em}P_{n}(\nu_{k}(P'))(
u_{1}\otimes \cdots \otimes u_{k}\otimes 
v ))).
\end{eqnarray}
Let $\tilde{\lambda}$ be an element of $V'$ defined by
$$\tilde{\lambda}(v)=\res_{w}(f_{A}(w))^{-1}f'_{A}(w)
\lambda(Y(e^{-\sum_{j>0}B(f_{A}(w))L(j)}
u_{0}, w)v).$$
Then by (\ref{gamma-k-2}) and (\ref{gamma-k-3}),
$\tilde{\lambda}\in \tilde{G}$ and thus 
the right-hand side of (\ref{gamma-k-1}) is in 
$F_{k}$.

\begin{prop}
The map $\gamma_{k}$ is continuous and satisfies 
\begin{equation}\label{halfinv}
\gamma_{k}
\circ \iota_{F_{k}}=I_{F_{k}}
\end{equation}
where $I_{F_{k}}$ is the identity map
on $F_{k}$. 
\end{prop}
\pf
We still use $$(z_{0}, \dots,
z_{k-1}; A, (1, \mathbf{0}), \dots, (1, \mathbf{0}))$$ instead of 
$$(z_{1}, \dots, z_{k}; A, (1, \mathbf{0}), \dots, (1, \mathbf{0}))$$ 
to denote an element of
$\Theta_{k+1}$. 
We know that there exists $t\in [0, 1)$ such that
for $\epsilon\in [t, 1]$, $\{w\in \mathbb{C}\;|\;|w|=\epsilon\}$
is in the domain of 
$f_{A}(1/w).$
Let
$C_{\epsilon}=f_{A}(\{w\in \mathbb{C}\;|\;|w|=1/\epsilon\})$ for 
$\epsilon\in [t, 1]$.
Then by the definition of $\gamma_{k}$ and Cauchy's theorem,
for any $\epsilon\in [t, 1]$ such that $z_{1}, \dots, z_{k-1}$
are in the interior of $C_{\epsilon}$, we have
\begin{eqnarray*}
\lefteqn{\gamma_{k}(g_{k+1}(\lambda\otimes
u_{0}\otimes u_{1}\otimes \cdots \otimes 
u_{l}\otimes v))}\nno\\
&&=
\frac{1}{2\pi \sqrt{-1}}\oint_{C_{\epsilon}}z_{0}^{-1}
g_{k+1}(\lambda\otimes u_{0}\otimes u_{1}\otimes \cdots \otimes 
u_{k}\otimes v)dz_{0}
\end{eqnarray*}
for $\lambda\in \tilde{G}$, $u_{0}, \dots, u_{k}, 
v\in V$. Thus by the definition of $J_{n}^{(k)}$, for any $n>0$,
there exists $\epsilon_{n}\in [t, 1]$
such that 
\begin{eqnarray}\label{esti-1}
\lefteqn{\|\gamma_{k}(g_{k+1}(\lambda\otimes
 u_{0}\otimes u_{1}\otimes \cdots \otimes 
u_{k}\otimes v))\|_{F_{k}, n}}\nno\\
&&=\sup_{(z_{1}, \dots,
z_{k-1}; A, (1, \mathbf{0}), \dots, (1, \mathbf{0}))\in J_{n}^{(k)}}
|\gamma_{k}(g_{k+1}(\lambda\otimes u_{0}\otimes u_{1}\otimes \cdots \otimes 
u_{k}\otimes v))|\nno\\
&&=\sup_{(z_{1}, \dots,
z_{k-1}; A, (1, \mathbf{0}), \dots, (1, \mathbf{0}))\in J_{n}^{(k)}}\nno\\
&&\hspace{2em}
\biggl|\frac{1}{2\pi \sqrt{-1}}\oint_{z_{0}\in C_{\epsilon_{n}}}z_{0}^{-1}
g_{k+1}(\lambda\otimes u_{0}\otimes u_{1}\otimes \cdots \otimes 
u_{k}\otimes v)dz_{0}\biggr|\nno\\
&&\le \sup_{(z_{1}, \dots,
z_{k-1}; A, (1, \mathbf{0}), \dots, (1, \mathbf{0}))
\in J_{n}^{(k)}, z_{0}\in C_{\epsilon_{n}}}
|g_{k+1}(\lambda\otimes u_{0}\otimes u_{1}\otimes \cdots \otimes 
u_{k}\otimes v)|.\nno\\
&&
\end{eqnarray}
For any $z_{0}\in C_{\epsilon_{n}}$, it is clear that there always exists 
positive integer $n_{z_{0}}$ and a open subset $U_{z_{0}}$ of $\mathbb{C}$
containing $z_{0}$
such that $U_{z_{0}}\times J_{n}^{(k)}\subset J_{n_{z_{0}}}^{(k+1)}$.
Since $C_{\epsilon_{n}}$ is compact, there exists finitely many points 
$z_{0}^{(1)}, \dots, z_{0}^{(l)}\in C_{\epsilon_{n}}$ such that 
$U_{z_{0}^{(1)}}, \dots, U_{z_{0}^{(l)}}$ cover
$C_{\epsilon_{n}}$. 
Thus the right-hand side of (\ref{esti-1}) is less than or equal to
\begin{eqnarray}\label{esti-2}
&{\dps \sum_{i=1}^{l}\max_{(z_{0}, \dots,
z_{k-1}; A, (1, \mathbf{0}), \dots, (1, \mathbf{0}))\in J_{n_{z_{0}^{(i)}}}^{(k+1)}}
|g_{k+1}(\lambda\otimes u_{0}\otimes u_{1}\otimes \cdots \otimes 
u_{k}\otimes v)|}&\nno\\
&{\dps =\sum_{i=1}^{l}
\|g_{k+1}(\lambda\otimes  u_{0}\otimes u_{1}\otimes \cdots \otimes 
u_{k}\otimes v)\|_{F_{k+1}, n_{z_{0}^{(i)}}}.}&
\end{eqnarray}
Combining (\ref{esti-1}) and (\ref{esti-2}), we see that 
$\gamma_{k}$ is continuous.

For $\lambda\in \tilde{G}$, $u_{1}, \cdots,
u_{k}, v\in V$, by definition,
\begin{eqnarray*}
\lefteqn{g_{k+1}(\lambda\otimes  \mathbf{1}\otimes u_{1}\otimes \cdots \otimes 
u_{k}\otimes v)}\nno\\
&&=\iota_{F_{k}}(g_{k}(\lambda
  \otimes u_{1}\otimes \cdots \otimes 
u_{k}\otimes v)).
\end{eqnarray*}
By definition, 
$$g_{k+1}(\lambda\otimes 
\mathbf{1}\otimes u_{1}\otimes \cdots \otimes 
u_{k}\otimes v)=g_{k}(\lambda\otimes u_{1}\otimes \cdots \otimes 
u_{k}\otimes v).$$
Thus
\begin{eqnarray*}
\lefteqn{\gamma_{k}(g_{k+1}(\lambda\otimes \mathbf{1}
\otimes u_{1}\otimes \cdots \otimes 
u_{k}\otimes v))}\nno\\
&&=g_{k}(\lambda
  \otimes u_{1}\otimes \cdots \otimes 
u_{k}\otimes v).
\end{eqnarray*}
 So we have
(\ref{halfinv}).
\epfv

The proof of the following consequence is the same as the proof of 
Corollary 1.3 in \cite{H4}:

\begin{cor}\label{k*:k+1*}
The  adjoint 
map $\gamma_{k}^{*}$  of $\gamma_{k}$ 
satisfies 
\begin{equation}\label{halfinv*}
\iota_{F_{k}}^{*}\circ \gamma_{k}^{*}
=I_{F_{k}^{*}}
\end{equation}
where  
$$\iota_{F_{k}}^{*}: 
F_{k+1}^{*}\to F_{k}^{*}$$
is the adjoint  of $\iota_{F_{k}}$ and 
$I_{F_{k}^{*}}$ is the identity on $F_{k}^{*}$. It  
is
injective and continuous. As a map from $F_{k}^{*}$ to
$\gamma_{k}^{*}(F_{k}^{*})$, it is also open.
In particular, if we identify $F_{k}^{*}$ 
with $\gamma_{k}^{*}(F_{k}^{*})$,
the topology on $F_{k}^{*}$ is 
induced from the one on $F_{k+1}^{*}$.\epf
\end{cor}

In the rest of this section, we give the remaining steps in
the construction of the locally convex completion. 
These steps are mostly the same as those in \cite{H4}. 
Thus our description of these steps shall be brief.
Also we warn the reader again that although the notations we use
below are the same as those in \cite{H4}, they denote different things in the 
present paper.

We use $\langle \cdot, \cdot\rangle$ to denote the pairing between
$\tilde{G}$ and the algebraic dual space $\tilde{G}^{*}$ of 
$\tilde{G}$. It is an extension of the pairing between $V'$ and $V$
denoted using the same symbol. 
The spaces $\tilde{G}$ and $\tilde{G}^{*}$ with this pairing
form a dual pair 
of vector spaces and thus
give a locally convex topology to $\tilde{G}^{*}$.
The dual space $\tilde{G}^{*}$ can be viewed as a subspace of 
$(V')^{*}=\overline{V}$. 
We define 
$$e_{k}: V^{\otimes (k+1)}\otimes F_{k}^{*}\to
\tilde{G}^{*}\subset \overline{V}$$
by
$$
\langle \lambda, e_{k}(u_{1}\otimes \cdots \otimes u_{k}\otimes
v\otimes \mu)\rangle
=\mu(g_{k}(\lambda\otimes u_{1}\otimes \cdots \otimes u_{k}\otimes
v))
$$
for $\lambda\in \tilde{G}$, $u_{1}, \dots, u_{k}, v\in V$ and 
$\mu \in F_{k}^{*}$. 

We now have to assume that $V$ is finitely generated.
Let $X$ be the finite-dimensional subspace $X$ of $V$ 
spanned by a finite set of generators of $V$ containing the vacuum vector 
$\mathbf{1}$.  We give $X$ the topology induced by any norm on $X$.
Then
$X^{\otimes (k+1)}\otimes F^{*}_{k}$ is  a 
locally convex
space. Let $G_{k}$ be
the image $e_{k}(X^{\otimes (k+1)}\otimes F_{k}^{*})$ of 
$X^{\otimes (k+1)}\otimes F_{k}^{*}\subset V^{\otimes (k+1)}\otimes 
F_{k}^{*}$
under $e_{k}$. 

The proofs of Propositions \ref{inclusion} and \ref{elcont} below
are the same as the proofs of Propositions 1.4 and 1.5 in \cite{H4}:

\begin{prop}\label{inclusion}
For any $k\ge 0$, $G_{k}\subset G_{k+1}$.\epf
\end{prop}

\begin{prop}\label{elcont}
The linear map 
$$e_{k}|_{X^{\otimes (k+1)}\otimes F_{k}^{*}}: 
X^{\otimes (k+1)}\otimes F_{k}^{*}\to
\tilde{G}^{*}$$
is continuous. \epf
\end{prop}

\begin{cor}\label{quotient}
The quotient space 
$$(X^{\otimes (k+1)}\otimes
F^{*}_{k})/(e_{k}|_{X^{\otimes (k+1)}\otimes F^{*}_{k}})
^{-1}(0)$$
is a  locally convex
space. \epf
\end{cor}

Using the isomorphism from $G_{k}$ to 
$$(X^{\otimes (k+1)}\otimes
F_{k}^{*})/(e_{k}|_{X^{\otimes (k+1)}\otimes F_{k}^{*}})^{-1}(0),$$
we obtain a locally convex
space structure on $G_{k}$ from that   on 
$$(X^{\otimes (k+1)}\otimes
F_{k}^{*})/(e_{k}|_{X^{\otimes (k+1)}\otimes F_{k}^{*}})^{-1}(0).$$
Let $H_{k}$ be the completion of
$G_{k}$. Then $H_{k}$ is a complete locally convex space.

The proof of the following proposition is the same as the proof of 
Proposition 1.7 in \cite{H4}:

\begin{prop}\label{ind-tp}
The space $H_{k}$ can be embedded canonically in $H_{k+1}$. 
The topology on $H_{k}$ is the same as the one induced from the
topology on $H_{k+1}$. \epf
\end{prop}

Now we have a sequence $\{H_{k}\}_{k\ge 0}$ 
of strictly
increasing complete locally convex
spaces.
Let 
$$H=\bigcup_{k\ge 0}H_{k}$$
equipped with the inductive
limit topology. Then $H$ is a complete locally convex space. 
Let 
$$G=\bigcup_{k\ge 0}G_{k}\subset H.$$
Then $V\subset G$ and $G$ is  dense in $H$. The same argument as in 
\cite{H4} shows that  $G$ is in the
closure of $V$. Thus we have:

\begin{thm}
The vector space $H$ equipped with the strict inductive limit topology
is a locally convex 
completion of $V$.\epf
\end{thm}

\renewcommand{\theequation}{\thesection.\arabic{equation}}
\renewcommand{\thethm}{\thesection.\arabic{thm}}
\setcounter{equation}{0}
\setcounter{thm}{0}

\section{The locally convex completion and a semi-group of annuli}

In this section, we construct, on the topological
completion $H$, a structure of a representation of 
the semi-group of the $\frac{c}{2}$-th power   of the
 determinant line bundle  over the moduli space of 
conformal equivalence classes
of annuli with analytically parametrized 
boundary components.

Consider the moduli space $\mathcal{B}_{1,1,0}$
of annuli, that is,  the space of conformal 
equivalence classes of genus-zero Riemann surfaces with two boundary 
components, one positively oriented and one negatively oriented, and with 
analytic boundary parametrizations of the boundary components. 
There is a
sewing operation on $\mathcal{B}_{1,1,0}$ such that it becomes a semi-group.
(See Appendix D of \cite{H3} for details.) 
There is a determinant line bundle $\dt(1)$ 
over $\mathcal{B}_{1,1,0}$
and its $c$-th power $\dt^{c}(1)$ for any $c\in \mathbb{C}$
is well-defined.

\begin{prop}
For any complex number $c$,  $\dt^{c}(1)$  
 has a structure of a semi-group
and is the central extension of $\mathcal{B}_{1,1,0}$ with central 
charge $2c$. \epf
\end{prop}

This result and its proof are contained implicitly in Appendix D of \cite{H3}. 
See \cite{H3} for details. 

By the uniformization theorem, it is clear that the semi-group 
$\mathcal{B}_{1,1,0}$ is 
isomorphic to the semi-group of the moduli space $K_{\mathfrak{H}_{1}}(1)$
equipped with the sewing operation. We shall identify 
 $\mathcal{B}_{1,1,0}$
with $K_{\mathfrak{H}_{1}}(1)$. Over the moduli space $K(1)$, we have 
a determinant line bundle and its $\frac{c}{2}$-th power 
$\tilde{K}^{c}(1)$ for 
any complex number $c$. We denote the restriction 
of $\tilde{K}^{c}(1)$ to $K_{\mathfrak{H}_{1}}(1)$ by 
$\tilde{K}^{c}_{\mathfrak{H}_{1}}(1)$. Then $\tilde{K}^{c}_{\mathfrak{H}_{1}}(1)$
is a semi-group isomorphic to $\dt^{c/2}(1)$. 
See \cite{H3} for details. We now construct 
a structure of a representation of $\tilde{K}^{c}_{\mathfrak{H}_{1}}(1)$
on $H$ where $c$ is the central charge of $V$. 

First we give a right action of 
$\tilde{K}^{c}_{\mathfrak{H}_{1}}(1)$ on  
 $\tilde{G}$. Let $\lambda\in \tilde{G}$ and $\tilde{Q}=(Q; C)\in 
\tilde{K}^{c}_{\mathfrak{H}_{1}}(1)$ (where 
$Q\in K_{\mathfrak{H}_{1}}(1)$
and $C\in \mathbb{C}$). We define $\lambda_{\tilde{Q}}\in V^{*}$
by
\begin{equation}\label{action1}
\lambda_{\tilde{Q}}(v)=
C\sum_{n\in \mathbb{Z}}\lambda(P_{n}((\nu_{1}(Q))(v)).
\end{equation}
Note that the right-hand side of (\ref{action1}) is 
absolutely convergent because $\lambda\in \tilde{G}$. 

\begin{lemma}
The linear functional $\lambda_{\tilde{Q}}$ is in fact in 
$\tilde{G}$.
\end{lemma}
\pf
By definition, for any $P\in \Theta(k)$,
\begin{eqnarray}\label{action2}
\lefteqn{\sum_{n\in \mathbb{Z}}
\lambda_{\tilde{Q}}(P_{n}(Q(u_{1}, \dots, u_{k}, v; P)))}\nno\\
&&=\sum_{n\in \mathbb{Z}}
\lambda_{\tilde{Q}}(P_{n}((\nu_{k}(P))(u_{1}\otimes 
\cdots \otimes u_{k}\otimes 
v)))\nno\\
&&=C\sum_{n\in \mathbb{Z}}\sum_{m\in \mathbb{Z}}
\lambda(P_{m}((\nu_{1}(Q))(P_{n}((\nu_{k}(P))(u_{1}\otimes 
\cdots \otimes u_{k}\otimes 
v)))).
\end{eqnarray}
We want to show that the right-hand side of (\ref{action2}) 
is absolutely convergent. To show this convergence, we note that,
by the sewing axiom for geometric vertex operator algebras,
\begin{eqnarray}\label{action3}
\lefteqn{\sum_{m\in \mathbb{Z}}\sum_{n\in \mathbb{Z}}
\lambda(P_{m}((\nu_{1}(Q))(P_{n}((\nu_{k}(P))(u_{1}\otimes 
\cdots \otimes u_{k}\otimes 
v))))}\nno\\
&&=\sum_{m\in \mathbb{Z}}
\lambda(P_{m}(((\nu_{1}(Q))_{^{1}}*_{^{0}}(\nu_{k}(P)))(u_{1}\otimes 
\cdots \otimes u_{k}\otimes 
v))))\nno\\
&&=\sum_{m\in \mathbb{Z}}
\lambda(P_{m}((\nu_{k}(Q_{^{1}}\infty_{^{0}}P))(u_{1}\otimes 
\cdots \otimes u_{k}\otimes 
v))).
\end{eqnarray}
Note that since $\lambda\in \tilde{G}$, the right-hand side of (\ref{action3}) 
is absolutely convergent and is 
analytic in $P$ and $Q$. Thus the double sum and the iterated sum in the 
other order are also absolutely convergent. Since 
the iterated sum in the right-hand side of (\ref{action2}) 
is exactly the iterated sum in the 
other order, it is absolutely convergent. \epfv

By this lemma, $\lambda \mapsto \lambda_{\tilde{Q}}$ for 
$\lambda\in \tilde{G}$ give a right action 
of $\tilde{K}^{c}_{\mathfrak{H}_{1}}(1)$ on  
 $\tilde{G}$. This right action induces a left action on 
$\tilde{G}^{*}$. It also
 induces right actions of $\tilde{K}^{c}_{\mathfrak{H}_{1}}(1)$
on $F_{k}$, $k\ge 0$, as follows:
\begin{eqnarray*}
g_{k}(\lambda\otimes u_{1}\otimes \cdots \otimes u_{k}\otimes v)
&\mapsto& g^{\tilde{Q}}_{k}(\lambda \otimes u_{1}\otimes \cdots\otimes u_{k}
\otimes v)\nno\\
&&=g_{k}(\lambda_{\tilde{Q}} \otimes u_{1}\otimes \cdots\otimes u_{k}
\otimes v),
\end{eqnarray*}
for $\tilde{Q}\in \tilde{K}^{c}_{\mathfrak{H}_{1}}(1)$,
$\lambda\in \tilde{G}$, $u_{1}, \dots, u_{k}, v\in V$.
These right actions on $F_{k}$, $k\ge 0$, induce left actions on 
$F^{*}_{k}$. For simplicity,
we shall also use $\tilde{Q}$ to  denote the 
images of $\tilde{Q}\in \tilde{K}^{c}_{\mathfrak{H}_{1}}(1)$ in $\edo \tilde{G}^{*}$
and $\edo F^{*}_{k}$, $k\ge 0$.

\begin{prop}
For $k\ge 0$, $\tilde{Q}\in \tilde{K}^{c}_{\mathfrak{H}_{1}}(1)$,
$\mu\in F_{k}^{*}$, $u_{1}, \dots, u_{k}, v\in V$.
$$\tilde{Q}\cdot e_{k}(u_{1}\otimes \cdots \otimes u_{k}\otimes
v\otimes \mu)=e_{k}(u_{1}\otimes \cdots \otimes u_{k}\otimes
v\otimes \tilde{Q}\cdot \mu).$$
\end{prop}
\pf
This follows from the definitions of $e_{k}$ and the 
left actions of $\tilde{K}^{c}_{\mathfrak{H}_{1}}(1)$ on $\tilde{G}^{*}$
and $F^{*}_{k}$.
\epfv

By this proposition, we immediately obtain:

\begin{cor}\label{2-5}
For $k\ge 0$, the actions of $\tilde{K}^{c}_{\mathfrak{H}_{1}}(1)$ on $\tilde{G}^{*}$
and $F^{*}_{k}$ induce an action of $\tilde{K}^{c}_{\mathfrak{H}_{1}}(1)$
on $G_{k}$ and thus an action on $H_{k}$. The actions of 
$\tilde{K}^{c}_{\mathfrak{H}_{1}}(1)$
on $H_{k}$ induce an action on $H$.\epf
\end{cor}

We shall still use  $\tilde{Q}$ to denote the images of 
$\tilde{Q}\in \tilde{K}^{c}_{\mathfrak{H}_{1}}(1)$
in $\edo H_{k}$, $k\ge 0$, and $\edo H$. We have the following:

\begin{prop}\label{2-6}
Let $\tilde{P}\in \tilde{K}^{c}_{\mathfrak{H}_{1}}(1)$. Then 
its images in $\edo H_{k}$, $k\ge 0$, and $\edo H$ are continuous. 
\end{prop}
\pf
We need only  prove the continuity of  the images of $\tilde{Q}$
in $\edo H_{k}$, $k\ge 0$. Since the actions on $H_{k}$, $k\ge 0$, are 
induced from the action on $\tilde{G}^{*}$, we need only show that 
the image of $\tilde{Q}$ in $\edo \tilde{G}^{*}$ is continuous. 
This is equivalent to the continuity of the image of 
$\tilde{Q}$ in $\edo \tilde{G}$. But from the definition
(\ref{action1}), it is clear that 
the image of 
$\tilde{Q}$ in $\edo \tilde{G}$ is continuous. \epfv

Combining Corollary \ref{2-5} and Proposition \ref{2-6}, 
we obtain the following:

\begin{thm}
The complete locally convex spaces $H_{k}$, $k\ge 0$, and 
$H$ have structures of continuous representations of 
$\tilde{K}^{c}_{\mathfrak{H}_{1}}(1)$ or
of $\dt^{c/2}(1)$.\epf
\end{thm}

Note that in the constructions of $H_{0}$ and of the structure 
of a continuous representation of $\tilde{K}^{c}_{\mathfrak{H}_{1}}(1)$
on $H_{0}$, only the structure of a $\mathbb{Z}$-graded representation
of the Virasoro algebra on $V$ and a certain lower-truncation 
condition of the representation is used. Thus we actually have the following:

\begin{thm}
Let $V=\coprod_{n\in \mathbb{Z}}V_{(n)}$
be a $\mathbb{Z}$-graded module for the Virasoro 
algebra satisfying the  conditions: {\rm (i)} $L(0)v=nv$ for
$v\in V_{(n)}$ and {\rm (ii)} for any $v\in V$, the $\mathbb{Z}$-graded
submodule 
$W=\coprod_{n\in \mathbb{Z}}W_{(n)}$ for the 
Virasoro algebra generated by $v$ is lower truncated, that is,
$W_{(n)}=0$ when $n$ is sufficiently small. 
Then the same constructions in Section 1 and in this 
section gives a locally convex completion $H_{0}$ of $V$ and a
structure of a continuous 
representation of $\tilde{K}^{c}_{\mathfrak{H}_{1}}(1)$ or
of $\dt^{c/2}(1)$ on $H_{0}$.\epf
\end{thm}

\renewcommand{\theequation}{\thesection.\arabic{equation}}
\renewcommand{\thethm}{\thesection.\arabic{thm}}
\setcounter{equation}{0}
\setcounter{thm}{0}

\section{The locally convex completion and 
the vertex operator map}

Consider a conformal equivalence class of 
genus-zero Riemann surfaces with three ordered boundary components,
the first positively oriented and the other two negatively oriented,
and with analytic
parametrizations at these boundary components. Such a 
conformal equivalence class can be naturally identified with 
an element of $K_{\mathfrak{H}_{1}}(2)$ (see \cite{H3}). 
We shall denote the corresponding 
element in $K_{\mathfrak{H}_{1}}(2)$ by $Q$. 
Then a pair consisting of such a  
conformal equivalence class  and an element of the 
$\frac{c}{2}$-th power of
the determinant line over it corresponding to an element $\tilde{Q}$ of 
$\tilde{K}^{c}_{\mathfrak{H}_{1}}(2)$. 

In this section, we use the vertex operator map to
construct continuous linear maps from 
the topological completion of $H\otimes H$
to $H$ associated to $\tilde{Q}\in \tilde{K}^{c}_{\mathfrak{H}_{1}}(2)$.

Let $H\widetilde{\otimes} H$ be 
the locally convex completion of the vector space tensor product
$H\otimes H$. We would like to construct a continuous linear map 
$$\overline{\Psi}_{Y}(\tilde{Q}): H\widetilde{\otimes} H\to H$$ 
associated to $\tilde{P}$ such that restricting to 
$V\otimes V$, it is equal to the  linear map
$\Psi_{Y}(\tilde{Q}): 
V\otimes V
\to \overline{V}$
constructed in \cite{H3}. Because $\tilde{K}^{c}_{\mathfrak{H}_{1}}(2)$
is infinite-dimensional,  our construction here is 
more complicated than the one in \cite{H4}. Nevertheless,
the idea and the steps are
mostly the 
same. Because of this, we shall be brief in our arguments below.

Given any $Q\in K(2)$, let $Q'$ be the element of $K(2)$ 
obtained by switching the negatively oriented and the second positively 
oriented punctures of any sphere with tubes in $P$. 
Thus we obtain a bijective map $'$ from $K(2)$ to itself. 
Since the line bundle $\tilde{K}^{c}(2)$ is canonically trivial,
this map $'$ can be extended to a bijective map $'$ from 
$\tilde{K}^{c}(2)$ to itself. It is clear that this map $'$
map $\tilde{K}^{c}_{\mathfrak{H}_{1}}(2)$ to itself.

We now fix $\tilde{Q}\in \tilde{K}^{c}_{\mathfrak{H}_{1}}(2)$.
For any $\lambda\in \tilde{G}$ and $u\in V$, we define an element 
$u\diamond_{\tilde{Q}}\lambda\in V^{*}$ by
$$
(u\diamond_{\tilde{Q}}\lambda)(v)
=\sum_{n\in \mathbb{Z}}
\lambda(P_{n}((\Psi_{2}(\tilde{Q}'))(u\otimes v))
$$
for $v\in V$.

\begin{prop}\label{3-1}
The element $u\diamond_{\tilde{Q}}\lambda$ is in $\tilde{G}$.
\end{prop}
\pf
We write 
$\tilde{Q}=(Q; C)$. For any $k\ge 0$, $u_{1}, \dots, u_{k},
v\in V$, $P\in \Theta_{k}$,
\begin{eqnarray}\label{3.1}
\lefteqn{\sum_{m\in \mathbb{Z}} 
(u\diamond_{\tilde{Q}}\lambda)(P_{m}(Q(u_{1}, \dots, u_{k}, v;
P)))}\nno\\ 
&&=\sum_{m\in \mathbb{Z}}\sum_{n\in
\mathbb{Z}} \lambda(P_{n}((\Psi_{2}(\tilde{Q}'))(u\otimes
P_{m}(Q(u_{1}, \dots, u_{k}, v;
P))))\nno\\
&&=C\sum_{m\in \mathbb{Z}}\sum_{n\in
\mathbb{Z}} \lambda(P_{n}((\nu_{2}(Q'))(u\otimes 
P_{m}(Q(u_{1}, \dots, u_{k}, v;
P))))\nno\\
&&=C\sum_{m\in \mathbb{Z}}\sum_{n\in
\mathbb{Z}} \lambda(P_{n}((\nu_{2}(Q'))(u\otimes 
P_{m}(Q(u_{1}, \dots, u_{k}, v;
P))))\nno\\
&&=C\sum_{m\in \mathbb{Z}}\sum_{n\in
\mathbb{Z}} \lambda(P_{n}(Q(u,
P_{m}(Q(u_{1}, \dots, u_{k}, v;
P)); Q'))).
\end{eqnarray}
We need to prove that the right-hand side of (\ref{3.1}) is 
absolutely convergent.

As in \cite{H4}, we consider the iterated sum in the other order
$$C\sum_{n\in \mathbb{Z}}\sum_{m\in
\mathbb{Z}} \lambda(P_{n}(Q(u,
P_{m}(Q(u_{1}, \dots, u_{k}, v;
P)); Q')))$$
which is convergent by using the sewing axiom and 
the fact that $\lambda\in \tilde{G}$. Moreover
it is clear that this iterated sum is the expansion 
of an analytic function in two variables evaluated 
at a certain particular point. Thus the double sum 
must be absolutely convergent and consequently the
right-hand side of (\ref{3.1}) is absolutely convergent.
\epfv

For any $l\ge 0$, 
we define a linear map 
$\alpha_{l}: \tilde{G} \otimes X^{l+1}\otimes F_{l}^{*}\to V^{*}$ by
\begin{eqnarray*}
\lefteqn{(\alpha_{l}(\lambda\otimes v_{1}\otimes \dots 
\otimes v_{l}\otimes v\otimes \mu))(u)}\nno\\
&&=\langle u\diamond_{\tilde{Q}}\lambda,
e_{l}
(v_{1}\otimes \dots 
\otimes v_{l}\otimes v\otimes \mu)\rangle.
\end{eqnarray*}
for $\lambda\in \tilde{G}$,  
$v_{1}, \dots, v_{l}, v\in X$,
$\mu\in F_{l}^{*}$ and $u\in V$.

\begin{prop}\label{3-2}
The image of  $\alpha_{l}$ is  in $\tilde{G}$.
\end{prop}
\pf
For any $k\ge 0$,
$\lambda\in \tilde{G}$, $u_{1}, \dots, u_{k}, u\in V$, $P\in \Theta_{k}$,
$v_{1}, \dots, v_{l},  v\in X$
and 
$\mu\in F_{l}^{*}$,
\begin{eqnarray}\label{3.2}
\lefteqn{\sum_{n\in \mathbb{Z}}
(\alpha_{l}(\lambda\otimes v_{1}\otimes \dots 
\otimes v_{l}\otimes v\otimes \mu))(P_{n}(Q(u_{1}, \dots, u_{k}, u; 
P)))}\nno\\
&&=\sum_{n\in \mathbb{Z}}\langle (P_{n}(Q(u_{1}, \dots, u_{k}, u; 
P))\diamond_{\tilde{Q}}\lambda), e_{l}
(v_{1}\otimes \cdots \otimes v_{l} \otimes
v\otimes \mu)\rangle\nno\\
&&=\sum_{n\in \mathbb{Z}}\mu(g_{l}((P_{n}(Q(u_{1}, \dots, u_{k}, u; 
P))\diamond_{\tilde{Q}}\lambda)\otimes 
v_{1}\otimes
\cdots \otimes \otimes v_{l}\otimes v))\nno\\
&&=\sum_{n\in \mathbb{Z}}
\mu\Biggl(\sum_{m\in \mathbb{Z}}(P_{n}(Q(u_{1}, \dots, u_{k}, u; 
P))\diamond_{\tilde{Q}}\lambda)(
P_{m}(Q(v_{1},
\dots, v_{k},
v; \cdot)))\Biggr)\nno\\
&&=\sum_{n\in \mathbb{Z}}\sum_{m\in \mathbb{Z}}\sum_{p\in \mathbb{Z}}
\mu(\lambda(P_{p}((\Psi_{2}(\tilde{Q}'))\nno\\
&&\hspace{6em}(
(P_{n}(Q(u_{1}, \dots, u_{k}, u; 
P))\otimes P_{m}(Q(v_{1},
\dots, v_{k},
v; \cdot))))\nno\\
&&=C\sum_{n\in \mathbb{Z}}\sum_{m\in \mathbb{Z}}\sum_{p\in \mathbb{Z}}
\mu(\lambda(P_{p}((\nu_{2}(\tilde{Q}'))\nno\\
&&\hspace{6em}((P_{n}(Q(u_{1}, \dots, u_{k}, u; 
P))\otimes P_{m}(Q(v_{1},
\dots, v_{k},
v; \cdot)))).\nno\\
&&
\end{eqnarray}
We need only to show that the right hand side of (\ref{3.2}) 
is absolutely convergent. The proof is similar to the 
proof in Proposition \ref{3-1} above: We first show that 
one of the iterated sums in other orders is absolutely convergent 
and is convergent to an analytic function in 
$Q$. Then this function can be expanded as series and the 
series is triply absolutely convergent. In particular, the 
iterated sum  in the right-hand side of (\ref{3.2}) 
is absolutely convergent and is equal to this triple 
sum. 
\epfv

By Proposition \ref{3-2}, 
$$\sum_{n\in \mathbb{Z}}\alpha_{l}(\lambda \otimes v_{1}\otimes \dots 
\otimes v_{l}\otimes v\otimes \mu)(P_{n}(Q(u_{1}, \dots, u_{k}, u; 
\tilde{Q})))$$
is absolutely convergent and  equal to 
$$g_{k}(\alpha_{l}(\lambda\otimes v_{1}\otimes \dots 
\otimes v_{l}\otimes v\otimes \mu)\otimes u_{1}\otimes \cdots\otimes 
u_{k}\otimes u)\in F_{k}.$$
We define a linear map 
$$\beta_{k, l}: F_{k}^{*}\otimes F_{l}^{*}
\to F_{k+l+1}^{*}$$ 
by
\begin{eqnarray*}
\lefteqn{(\beta_{k, l}(\mu_{1}, \mu_{2}))(g_{k+l+1}(\lambda
\otimes 
u_{1}\otimes \dots \otimes u_{k+1}\otimes u\otimes 
v_{1}\otimes \cdots\otimes v_{l}\otimes v))}\nno\\
&&=\mu_{1}
(g_{k}(\alpha_{l}(\lambda \otimes v_{1}\otimes \dots 
\otimes v_{l}\otimes v\otimes \mu_{2})\otimes u_{1}\otimes \cdots\otimes 
u_{k}\otimes u))
\end{eqnarray*}
for $\lambda\in \tilde{G}$, $u_{1}, \dots,  u_{k}, u, v_{1}, \dots,
v_{l}, v\in V$, $\mu_{1}\in F_{k}^{*}$ and $\mu_{2}\in F_{l}^{*}$. In fact
this formula only gives a linear map from $F_{k}^{*}\otimes
F_{l}^{*}$ to the algebraic dual of $F_{k+l+1}$.
The proof of the following result is completely analogous to Proposition 
2.3 in \cite{H4}:

\begin{prop}\label{3-4}
The image of the map $\beta_{k, l}$  is indeed in $F_{k+l+1}^{*}$ and the map
$\beta_{k, l}$ is continuous.\epf
\end{prop}

Let 
$$h_{1}=e_{k}(u_{1}\otimes \cdots \otimes u_{k}\otimes
u\otimes \mu_{1})\in G_{k}$$
and 
$$h_{2}=e_{l}(v_{1}\otimes \cdots \otimes v_{l}\otimes
v\otimes \mu_{2})\in G_{l}$$
where $u_{1}, \dots, u_{k}, u, v_{1}, \dots, 
v_{l}, v\in X$, 
$\mu_{1}\in F_{k}^{*}$ and $\mu_{2}\in 
F_{l}^{*}$.
We define 
\begin{eqnarray*}
\lefteqn{(\overline{\Psi}_{Y}(\tilde{Q}))(h_{1}\otimes h_{2})}\nno\\
&&=e_{k+l+1}(u_{1}\otimes \cdots \otimes 
u_{k}\otimes u\otimes v_{1}\otimes \cdots \otimes 
v_{l}\otimes  v\otimes \beta_{k, l}(\mu_{1}, \mu_{2})).
\end{eqnarray*}
Note that any element of $G_{k}$ or $G_{l}$ is a linear
combination of elements of the form
$h_{1}$ or $h_{2}$, respectively,  given above, and that $k$ and $l$ are
arbitrary. Thus we obtain a linear map 
$$\overline{\Psi}_{Y}(\tilde{Q})\mbar_{G\otimes G}: G\otimes G\to G.$$
The proof of the following result is completely analogous to the proof
of Proposition 
2.4 in \cite{H4}:

\begin{prop}\label{contong}
The map $\overline{\Psi}_{Y}(\tilde{Q})\mbar_{G\otimes G}$ is continuous.\epf
\end{prop}

Since $G$ is dense in $H$, we can extend 
$\overline{\Psi}_{Y}(\tilde{Q})\mbar_{G\otimes G}$ to a linear map 
$\overline{\Psi}_{Y}(\tilde{Q})$ from $H\widetilde{\otimes} H$ to $H$. 
The proof of the following theorem is completely analogous to the 
proof of Theorem2.5 
in \cite{H4}:

\begin{thm}
The map $\overline{\Psi}_{Y}(\tilde{Q})$ is a continuous extension of 
$\Psi_{Y}(\tilde{Q})$
to $H\widetilde{\otimes}H.$
That is, $\overline{\Psi}_{Y}(\tilde{Q})$ is continuous and
$$\overline{\Psi}_{Y}(\tilde{Q})\mbar_{V\otimes
V}=\Psi_{Y}(\tilde{Q}).\epfe$$
\end{thm}

\renewcommand{\theequation}{\thesection.\arabic{equation}}
\renewcommand{\thethm}{\thesection.\arabic{thm}}
\setcounter{equation}{0}
\setcounter{thm}{0}

\section{Locally convex  completions, operads and double loop spaces}

In this section, we reformulate the result  obtained in \cite{H4} and in
Sections 2 and 3 above using the language of operads.

First, the result in Section 2 of \cite{H4} immediately gives the following:

\begin{thm}\label{part1-main}
Let $V$ be a finitely-generated $\mathbb{Z}$-graded vertex
algebra.  Then the topological completion $H$ of $V$ constructed in
\cite{H4} has a structure of an algebra over the framed little disk
operad such that for the unit disk with two embedded disks of radius
$r_{1}$ and $r_{2}$ centered at $0$ and $z$, respectively, the
corresponding map from $H\widetilde{\otimes}H$ to $H$ is the map
$\overline{\nu}_{Y}([D(z, r_{1}, r_{2})])$. (See \cite{H4} for the notation
$\overline{\nu}_{Y}$ and $[D(z, r_{1}, r_{2})]$.)
\end{thm} 
\pf
The framed little disk operad is generated by 
the unit disk with two embedded disks of radius
$r_{1}$ and $r_{2}$ centered at $0$ and $z$ and 
the unit disk with the unit disk itself embedded and with 
the frames given by complex numbers $a$ of absolute value equal to $1$.
So we need only define the maps corresponding to these 
elements of the operad. For
the unit disk with two embedded disks of radius
$r_{1}$ and $r_{2}$ centered at $0$ and $z$, we define the associated map 
to be $\nu_{Y}([D(z, r_{1}, r_{2})])$. For
the unit disk with the unit disk itself embedded and with 
the frames given by complex numbers $a$ of absolute value equal to $1$,
we define the associated map 
to be $a^{L(0)}: H\to H$. Then we get a structure of algebra on $H$ 
over 
the framed little disk operad.
\epfv

Next, combining the results of \cite{H3} and 
the results in Sections 2 and 3 above, we obtain the following 
result:

\begin{thm}\label{main}
Let $V$ be a finitely-generated  vertex operator
algebra.  Then the topological completion $H$ of $V$ constructed 
above has a structure of an algebra over the operad 
$\tilde{K}^{c}_{\mathfrak{H}_{1}}$ or, 
equivalently, 
of $\dt^{c/2}$.  \epf
\end{thm} 

\begin{cor}\label{spc-l-d-o}
Let $V$ be a finitely-generated $\mathbb{Z}$-graded vertex 
algebra or a finitely-generated vertex operator algebra. Then 
locally convex completion $H$ of $V$ constructed in Part I (\cite{H4})
or in Section 1 above has a structure of a space over the 
little framed disk operad. In particular, it has 
a structure of a space over the 
little framed disk operad.
\end{cor}
\pf 
Since we have a natural continuous map from 
$H\times H$ to $H\otimes H$, we see from Theorem \ref{part1-main}
that when $V$ is a finitely-generated $\mathbb{Z}$-graded vertex 
algebra, its locally convex completion constructed in Part I has 
a structure of a space over the little framed disk operad.

If $V$ is a finitely-generated vertex operator algebra. 
Then note that the little framed disk operad can in fact be viewed as a 
suboperad of $K_{\mathfrak{H}_{1}}$. Also note that
 the sewing of the determinant lines over elements in the 
little framed disk operad is trivial (see Appendix D of \cite{H3}).
Thus $H$ has a structure of an algebra over the framed little disk operad
and consequently has a structure of a space over the 
little framed disk operad. \epfv

A subspace of a Hausdorff 
space is said to be {\it compactly closed}
if the intersection of the subspace with each compact subset 
of the Hausdorff space is closed. A Hausdorff space is said to be
{\it compactly generated} if every compactly closed subspace 
is closed. See \cite{St} (and 
\cite{W} and \cite{M2}) for the notion of compactly 
generated topological space and properties of these spaces. 
In \cite{M}, May proved, among other things, the following recognition
principle for double loop spaces:

\begin{thm}\label{may}
If a compactly generated Hausdorff based topological space 
has a structure of a space over the little disk operad,
then it has the weak homotopy of a double loop space.\epf
\end{thm}

From \cite{St} (see also \cite{W} and \cite{M2}), we know that 
we can make a Hausdorff space into a compactly generated 
Hausdorff space by giving it a new topology in which a subspace is 
closed if and only if it is compactly closed in the original 
topology. Since this functor is usually denoted by $k$,
here we call the space with the new compactly generated topology the 
{\it $k$-ification} of the original space. Note that 
in the category of compactly generated spaces, the product 
of spaces is defined to be the $k$-ification of the usual product
(see \cite{St}, \cite{W} and \cite{M2}).
The following lemma follows immediately from the properties of
$k$-ifications of topological spaces:

\begin{lemma}\label{k-ify-sp}
If a Hausdorff based topological space is a space over 
the little disk operad (with the usual products of topological spaces),
then the $k$-ification of the space has a natural structure of 
a space over the little disk operad (with the products of 
compactly generated spaces).\epf
\end{lemma}

Combining Corollary \ref{spc-l-d-o}
with Theorem \ref{may} and Lemma \ref{k-ify-sp}, we obtain:

\begin{thm}
The  $k$-ifications of the
locally convex completions constructed in \cite{H4}
and in Section 1 above have weak homotopy types of double loop spaces.\epf
\end{thm}

\renewcommand{\theequation}{\thesection.\arabic{equation}}
\renewcommand{\thethm}{\thesection.\arabic{thm}}
\setcounter{equation}{0}
\setcounter{thm}{0}

\section{Locally-grading restricted conformal 
vertex algebras and topological completions}

The results in the present paper are true also for algebras 
which do not satisfy the (global) grading-restriction conditions.  
We first need the following:

\begin{defn}
{\rm A {\it conformal vertex algebra}  is a $\mathbb{Z}$-graded 
vertex algebra equipped with a Virasoro element $\omega$
satisfying all the axioms for vertex operator algebras
except the two grading-restriction axioms. A conformal vertex 
algebra is said to be {\it locally grading-restricted} if 
for any element of the conformal vertex
algebra, the module $W=\coprod_{n\in \mathbb{Z}}W_{(n)}$ for the
Virasoro algebra generated by this element satisfies the
grading-restriction conditions, that is, $\dim W_{(n)}<\infty$ for $n\in
\mathbb{Z}$ and $W_{(n)}=0$ for $n$ sufficiently small. }
\end{defn}

\begin{rema}
{\rm In fact,  it is not difficult to show that the 
condition $\dim W_{(n)}<\infty$ in the definition 
above can be derived as a consequence. Thus for concrete examples,
one need only verify the lower-truncation condition 
$W_{(n)}=0$ for $n$ sufficiently small.}
\end{rema}

We have the following:

\begin{thm}\label{cva}
The constructions and results in \cite{H4} and in Sections 1, 2, 3 and 4 above
hold for finitely-generated locally-grading-restricted
conformal vertex algebras.
\end{thm}
\pf
Note that the constructions and results in \cite{H4} and in 
Sections 1, 2, 3 and 4 above need only the locally grading-restriction 
conditions: All the properties of vertex operator algebras used,
for example, commutativity, associativity,
rationality and the factorization of exponentials of infinite sums
of Virasoro operators, still hold if  the locally-grading-restriction 
conditions are satisfied. The details are left to 
the reader as an  exercise.
\epf

\begin{rema}
{\rm Theorem \ref{cva} has  been used in \cite{HZ}.}
\end{rema}

\renewcommand{\theequation}{\thesection.\arabic{equation}}
\renewcommand{\thethm}{\thesection.\arabic{thm}}
\setcounter{equation}{0}
\setcounter{thm}{0}

\section{A locally convex completion of a finitely-gen\-erated module 
and operads}

We give the results for modules in this section. 
Since the constructions and 
proofs are all similar to the case of algebras, we shall only 
state the final results. All the constructions and proofs are left to the 
reader as exercises.

\begin{thm}
Let $V$ be a finitely-generated 
vertex operator algebra of central charge $c$, 
$H$ its locally convex topological 
completion constructed in Section 1  and $W$ a 
finitely-generated $V$-module.  Then constructions completely 
analogous to those in Sections 1, 2, 3 and 4 above give a locally convex 
topological completion $H^{W}$ of $W$ and  a structure of 
a module for the algebra $H$ over the operad 
$\tilde{K}^{c}_{\mathfrak{H}_{1}}$ (or
equivalently
of $\dt^{c/2}$) on $H^{W}$. \epf
\end{thm}

\noindent {\small \sc Department of Mathematics, Kerchof Hall,
University of Virginia, Charlottesville, VA 22904-4137}
\vskip .5em
\noindent {\it and}
\vskip .5em
\noindent {\small \sc Department of Mathematics, Rutgers University,
110 Frelinghuysen Rd., Piscataway, NJ 08854-8019 (permanent address)}

\vskip 1em

\noindent {\em E-mail address}: yzhuang@math.rutgers.edu

\end{document}